\documentstyle[twoside]{article}
\title{A note on a hypergeometric transformation formula due to Slater with an application}
\author{
Y. S. Kim,\footnote{Department of Mathematics Education, Wonkwang University, Iksan, Korea
E-mail: yspkim@wonkwang.ac.kr}
\ \ Arjun. K. Rathie\footnote{School of Mathematical and Physical Sciences, Central University of Kerala,
Periye P.O. Dist. Kasaragad 671123, Kerala State, India.
E-Mail: akrathie@cukerala.edu.in} \ \ and 
\ R. B. Paris\footnote{School of Computing, Engineering and Applied Mathematics, University of Abertay Dundee, Dundee DD1 1HG, UK.
E-Mail: r.paris@abertay.ac.uk}\ \footnote{Corresponding author}
 \\}
\begin{document}
\def\f#1#2{\mbox{${\textstyle \frac{#1}{#2}}$}}
\def\dfrac#1#2{\displaystyle{\frac{#1}{#2}}}
\def\boldal{\mbox{\boldmath $\alpha$}}
\newcommand{\bee}{\begin{equation}}
\newcommand{\ee}{\end{equation}}
\newcommand{\lam}{\lambda}
\newcommand{\ka}{\kappa}
\newcommand{\al}{\alpha}
\newcommand{\th}{\theta}
\newcommand{\om}{\omega}
\newcommand{\Om}{\Omega}
\newcommand{\fr}{\frac{1}{2}}
\newcommand{\fs}{\f{1}{2}}
\newcommand{\g}{\Gamma}
\newcommand{\br}{\biggr}
\newcommand{\bl}{\biggl}
\newcommand{\ra}{\rightarrow}
\newcommand{\mbint}{\frac{1}{2\pi i}\int_{c-\infty i}^{c+\infty i}}
\newcommand{\mbcint}{\frac{1}{2\pi i}\int_C}
\newcommand{\mboint}{\frac{1}{2\pi i}\int_{-\infty i}^{\infty i}}
\newcommand{\gtwid}{\raisebox{-.8ex}{\mbox{$\stackrel{\textstyle >}{\sim}$}}}
\newcommand{\ltwid}{\raisebox{-.8ex}{\mbox{$\stackrel{\textstyle <}{\sim}$}}}
\renewcommand{\topfraction}{0.9}
\renewcommand{\bottomfraction}{0.9}
\renewcommand{\textfraction}{0.05}
\newcommand{\mcol}{\multicolumn}
\date{}
\maketitle
\begin{abstract}
In this note we state (with minor corrections) and give an alternative proof of a very general hypergeometric transformation formula due to Slater. 
As an application, we obtain a new hypergeometric transformation formula for a ${}_5F_4(-1)$ series with one pair of parameters differing by unity expressed as a linear combination of two ${}_3F_2(1)$ series.

\vspace{0.4cm}

\vspace{0.3cm}

\end{abstract}

\vspace{0.3cm}

\noindent{\bf 1. \  Introduction}
\setcounter{section}{1}
\setcounter{equation}{0}
\renewcommand{\theequation}{\arabic{section}.\arabic{equation}}
\vspace{0.3cm}

The generalized hypergeometric function with $p$ numeratorial and $q$ denominatorial parameters is defined by the series \cite[p.~41]{S}
\bee\label{e1}
{}_pF_q\left[\!\!\begin{array}{c} a_1, a_2, \ldots ,a_p\\b_1, b_2, \ldots ,b_q\end{array}\!; z\right]=
\sum_{n=0}^\infty \frac{(a_1)_n (a_2)_n \ldots (a_p)_n}{(b_1)_n (b_2)_n \ldots (b_q)_n}\,\frac{z^n}{n!},
\ee
where for nonnegative integer $n$ the Pochhammer symbol (or ascending factorial) is defined by $(a)_0=1$ and for $n\geq 1$ by $(a)_n=a(a+1)\ldots (a+n-1)$. 
When $p\leq q$ the above series on the right-hand side of (\ref{e1}) converges for $|z|<\infty$, but when $p=q+1$ convergence occurs when $|z|<1$ (unless the series terminates).

By employing Bailey's transform of double series, Slater \cite[p.~60, Eq. (2.4.10)]{S} derived the following very general hypergeometric transformation formula (written here in corrected form)
\[\sum_{n=0}^\infty \frac{((a))_n ((d))_n ((v))_{2n}}{((h))_n ((g))_n ((f))_{2n} \,n!}x^ny^nz^{2n} \times
{}_{U+D+V}F_{E+F+G}\left[\begin{array}{c} (u), (d)+n, (v)+2n\\(e), (g)+n, (f)+2n\end{array}\!;xwz\right]\]
\[=\sum_{n=0}^\infty \frac{((d))_n ((u))_n ((v))_n}{((e))_n ((f))_n ((g))_n\, n!}x^nw^nz^n \hspace{8cm}\]
\bee\label{e2}
\times
{}_{A+E+V+1}F_{U+H+F}\left[\begin{array}{c}-n, (a), 1-n-(e), (v)+n\\(h), 1-n-(u), (f)+n\end{array}\!;(-1)^{1+E-U} w^{-1}yz\right],
\ee
where we have adopted the convention of writing the finite sequence of parameters $(a_1, \ldots ,a_A)$ simply by $(a)$ and the product of Pochhammer symbols by
$$((a))_n\equiv (a_1)_n\ldots (a_A)_n,$$
where an empty product is understood to be unity. 
The general result (\ref{e2}) contains as special cases very many relationships between generalized hypergeometric functions.

Several interesting special cases of (\ref{e2}) have been obtained by Exton \cite{E2,E1}. In particular,
he gave the transformation \cite{E1}
\[\sum_{n=0}^\infty \frac{((g))_n (c)_n (d)_n}{((h))_n (f)_{2n}}\frac{x^ny^n}{n!} \times
{}_{2}F_{1}\left[\begin{array}{c} c+n, d+n\\f+2n\end{array}\!;x\right]\hspace{4cm}\]
\bee\label{e3}
\hspace{4cm}=\sum_{n=0}^\infty \frac{(c)_n (d)_n}{(f)_n}\,\frac{x^n}{n!}
\times
{}_{G+1}F_{H+1}\left[\begin{array}{c}-n,\\f+n,\end{array}\begin{array}{c}(g)\\ (h)\end{array}\!;-y\right].
\ee
This result follows from (\ref{e2}) by setting the parameters $(e)$, $(g)$, $(u)$, $(v)=0$, $(d)=(c,d)$, replacing $(a)$ by $(g)$ and letting $w=z=1$.

Hypergeometric identities and transformation formulas have wide applications in numerous areas of mathematics including in series systems of symbolic computer algebra manipulation. A list of such useful identities and transformation formulas can be found in Slater's book \cite[Appendix III]{S}.
We have the following known hypergeometric identities:
\bee\label{e11}
{}_4F_3\left[\begin{array}{c}a, 1+\fs a, b, c\\\fs a, 1+a-b, 1+a-c\end{array};1\right]=
\frac{\g(1+a-b) \g(1+a-c) \g(\fs a+\fs) \g(\fs a-b-c+\fs)}{\g(1+a) \g(1+a-b-c) \g(\fs a-b+\fs) \g(\fs a-c+\fs)}
\ee
provided $\Re (a-2b-2c)>-1$, and
\[{}_5F_4\left[\begin{array}{c}a, 1+\fs a, b, c, d\\\fs a, 1+a-b, 1+a-c, 1+a-d\end{array};1\right]\hspace{5cm}\]
\bee\label{e12}\hspace{4cm}=
\frac{\g(1+a-b) \g(1+a-c) \g(1+a-d) \g(1+a-b-c-d)}{\g(1+a) \g(1+a-b-c) \g(1+a-b-d) \g(1+a-c-d)}
\ee provided $\Re (a-b-c-d)>-1$. If we put $c=-n$ in (\ref{e11}) and $d=-n$ in (\ref{e12}), where $n$ is a non-negative integer, we have the terminating forms
\bee\label{e13}
{}_4F_3\left[\begin{array}{c}-n, a, 1+\fs a, b \\\fs a, 1+a-b, 1+a-c\end{array};1\right]=
\frac{(1+a)_n (\fs+\fs a-b)_n}{(\fs+\fs a)_n (1+a-b)_n}
\ee
and
\bee\label{e14}
{}_5F_4\left[\begin{array}{c}-n, a, 1+\fs a, b, c\\\fs a, 1+a-b, 1+a-c, 1+a-d\end{array};1\right]
=\frac{(1+a)_n (1+a-b-c)_n}{(1+a-b)_n (1+a-c)_n}.
\ee

Furthermore, by taking $x=1$ in (\ref{e3}) and making use of Gauss' summation theorem, 
Exton \cite{E1} obtained the following transformation formula
\[{}_{G+2}F_{H+2}\left[\begin{array}{c}(g),\\(h),\end{array}
\begin{array}{c} c, d\\f-c, f-d\end{array};y\right]\hspace{5cm}\]
\bee\label{e15}
=\frac{\g(f-c) \g(f-d)}{\g(f) \g(f-c-d)} \sum_{n=0}^\infty \frac{(c)_n (d)_n}{(f)_n n!}\,{}_{G+1}F_{H+1}\left[\begin{array}{c}-n,\\f+n,\end{array}\begin{array}{c}(g)\\(h)\end{array}; -y\right]
\ee
When $G=H+1$, this result holds for $|y|\leq 1$ when the parameters are such that the series on the left is defined and is convergent.
Using the results (\ref{e13}) and (\ref{e14}) in (\ref{e15}), he then deduced the following known hypergeometric identity:
\[{}_5F_4\left[\begin{array}{c}a, 1+\fs a, b, c, d\\ \fs a, 1+a-b, 1+a-c, 1+a-d\end{array}; -1\right]\hspace{4cm}\]
\bee\label{e16}
\hspace{3cm}=\frac{\g(1+a-c) \g(1+a-d)}{\g(1+a) \g(1+a-c-d)}\ {}_3F_2\left[\begin{array}{c}c, d, \fs+\fs a-b\\
\fs+\fs a, 1+a-b\end{array}; 1\right].
\ee
Again, it is tacitly assumed that the hypergeometric series in (\ref{e16}) is convergent. This requires that the parametric excess $s$, which is defined as the difference between the sum of denominator and numerator parameters, should satisfy $\Re (s)>-1$ when $y=-1$ and $\Re (s)>0$ when $y=1$.

Recently, some progress has been achieved in generalizing various hypergeometric identities; for this we refer to the papers cited in \cite{KRR, RR}. In 2010, Kim {\it et al.\/} \cite{KRR} generalized several well-known identities and in particular obtained the following result involving one pair of numeratorial and denominatorial parameters differing by unity:
\[{}_4F_3\left[\begin{array}{c}a, b, c, d+1\\1+a-b, 1+a-c, d\end{array};1\right]=\frac{\g(1+a-b) \g(1+a-c)}{\g(1+a) \g(1+a-b-c)}\]
\bee\label{e18}
\times \left\{\frac{a}{2d}\,\frac{\g(\fs+\fs a) \g(\fs+\fs a-b-c)}
{ \g(\fs+\fs a-b) \g(\fs+\fs a-c)}
+\left(1-\frac{a}{2d}\right)\,\frac{\g(1+\fs a) \ \g(1+\fs a-b-c)}
{\g(1+\fs a-b)\g(1+\fs a-c)}\right\}
\ee
provided $\Re (a-2b-2c)>-1$ and $d\neq 0, -1, -2, \ldots\ $. 
The result (\ref{e18}) may be regarded as a generalization of (\ref{e12}).

The aim in this note is to provide first another method of derivation of Slater's general transformation in (\ref{e2}) and to point out two significant misprints in \cite[p.~60, Eq. (2.4.10)]{S}.
As an application, we give
a generalization of the ${}_5F_4(-1)$ summation in (\ref{e16}) to the case when a pair of numeratorial and denominatorial parameters differ by unity. Our result will be established with the help of Exton's transformation formula (\ref{e15}) and the summation in (\ref{e18}).

\vspace{0.6cm}

\noindent{\bf 2. \ Derivation of (\ref{e2})}
\setcounter{section}{2}
\setcounter{equation}{0}
\renewcommand{\theequation}{\arabic{section}.\arabic{equation}}
\vspace{0.3cm}

\noindent In order to derive (\ref{e2}) we proceed as follows. Denoting the left-hand side of (\ref{e2}) by $S$ and
expressing the generalized hypergeometric functions  ${}_{U+D+V}F_{E+F+G}$ in its series form, we find after some simplification
\[S=\sum_{n=0}^\infty\sum_{m=0}^\infty\frac{((a))_n ((u))_m}{((h))_n ((e))_m}\,\frac{((d))_{n+m} ((v))_{2n+m}}{((g))_{n+m} ((f))_{2n+m}}\,\frac{x^{n+m}y^n w^mz^{2n+m}}{n!\, m!},\] 
where we have made use of the elementary identities
\[((\alpha))_n((\alpha)+n)_m=((\alpha))_{n+m},\qquad ((\alpha))_{2n}((\alpha)+2n)_m=((\alpha))_{2n+m}.\]
We now replace $m$ by $m-n$ and apply the result \cite[p.~56, Lemma 10(1)]{R}
\[\sum_{n=0}^\infty\sum_{k=0}^\infty A(k,n)=\sum_{n=0}^\infty\sum_{k=0}^n A(k,n-k)\]
for convergent double series, to obtain
\[S=\sum_{m=0}^\infty\sum_{n=0}^m\frac{((a))_n ((u))_{m-n}}{((h))_n ((e))_{m-n}}\,\frac{((d))_{m} ((v))_{n+m}}{((g))_{m} ((f))_{n+m}}\,\frac{x^{m}y^n w^{m-n}z^{n+m}}{n!\, (m-n)!}.\] 

Employing the identity $((\alpha))_{n+m}=((\alpha))_m((\alpha)+m)_n$, together with
\[(\alpha)_{m-n}=\frac{(-1)^n (\alpha)_m}{(1-\alpha-m)_n},\qquad (m-n)!=\frac{(-1)^n m!}{(-m)_n},\]
we then find
\[S=\sum_{m=0}^\infty\sum_{n=0}^m\frac{((d))_m ((v))_m ((v)+m)_n ((a))_n ((u))_{m-n}}{((g))_m((f))_m((f)+m)_n ((h))_n((e))_{m-n}}\,\frac{x^my^n w^{m-n}z^{m+n}}{n!\ (m-n)!}\]
\[=\sum_{m=0}^\infty\frac{((d))_m((v))_m((u))_m}{((g)_m((f))_m((e))_m}\,\frac{x^mw^mz^m}{m!}\hspace{5cm}\] \[\hspace{3cm}\times\sum_{n=0}^m\frac{(-m)_n(1-(e)-m)_n((a))_n((v)+m)_n}{(1-(u))_n((h))_n((f)+m)_n}(-1)^{(1+E-U)n}\frac{w^{-n}y^nz^n}{n!}\,.\]
Finally, expressing the inner series as a hypergeometric function, we easily arrive at the right-hand side of (\ref{e2}). This completes the proof of the transformation formula (\ref{e2}).\ \ \  $\Box$

\vspace{0.6cm}

\noindent{\bf 3. \ A new transformation formula}
\setcounter{section}{3}
\setcounter{equation}{0}
\renewcommand{\theequation}{\arabic{section}.\arabic{equation}}
\vspace{0.3cm}

\noindent The transformation formula for the ${}_5F_4(1)$ series with a pair of numeratorial and denominatorial parameters differing by unity to be established is given by the following theorem.
\newtheorem{theorem}{Theorem}
\begin{theorem}$\!\!\!.$\ \ For $e\neq 0, -1, -2, \ldots\,$,  the following summation holds true
\[{}_5F_4\left[\begin{array}{c}a,\ b,\ c,\ d,\ 1+e\\1+a-b, 1+a-c, 1+a-d, e\end{array};-1\right]\hspace{6cm}\]
\[\hspace{1cm}=\frac{\g(1+a-c) \g(1+a-d)}{\g(1+a) \g(1+a-c-d)}\left\{\frac{a}{2e}\,{}_3F_2\left[\begin{array}{c}
c, d, \fs+\fs a-b\\\fs+\fs a, 1+a-b\end{array};1\right] \right.\]
\bee\label{e21}
\hspace{4cm}\left.+\left(1-\frac{a}{2e}\right)\,{}_3F_2\left[\begin{array}{c} c, d, 1+\fs a-b\\1+\fs a, 1+a-b\end{array};1\right]\right\}
\ee
provided $\Re (a-c-d)>\max\{-1, \Re (b)-\f{3}{2}\}$.
\end{theorem}

\noindent{\it Proof.}\ \ \ Our derivation follows in a straightforward manner from Exton's transformation formula (\ref{e15}). If we take $y=-1$, $G=3$, $H=2$, $g_1=a$, $g_2=1+e$, $g_3=b$, $h_1=e$, $h_2=1+a-b$ and $f=1+a$ in  (\ref{e15}), we obtain after some simplification
\[F\equiv {}_5F_4\left[\begin{array}{c}a,\ b,\ c,\ d,\ 1+e\\1+a-b, 1+a-c, 1+a-d, e\end{array};-1\right]\hspace{4cm}\]
\bee\label{e22}
\hspace{1cm}=\frac{\g(1+a-c) \g(1+a-d)}{\g(1+a) \g(1+a-c-d)}\sum_{n=0}^\infty\frac{(c)_n (d)_n}{(1+a)_n n!}\,{}_4F_3\left[\begin{array}{c}-n, a, b, 1+e\\1+a-b, 1+a+n, e\end{array};1\right].
\ee

The terminating ${}_4F_3(1)$ series appearing on the right-hand side of (\ref{e22}) can be evaluated with the help of the summation in (\ref{e18}).
If we take $c=-n$ in this latter summation formula, where $n$ is a non-negative integer, we have
\[{}_4F_3\left[\begin{array}{c}-n, a, b, d+1\\1+a-b, 1+a+n, d\end{array};1\right]\hspace{6cm}\]
\bee\label{e19}
\hspace{1cm}=\frac{(1+a)_n}{(1+a-b)_n}\left\{\frac{a}{2d}\,\frac{(\fs+\fs a-b)_n}{(\fs+\fs a)_n}+\left(1-\frac{a}{2d}\right)\,\frac{(1+\fs a-b)_n}{(1+\fs a)_n}\right\}.
\ee
It is of interest to mention parenthetically that, when $d=\fs a$, (\ref{e18}) and (\ref{e19}) reduce to (\ref{e11}) and (\ref{e13}) respectively. 
Then we find that
\[F=\frac{\g(1+a-c) \g(1+a-d)}{\g(1+a) \g(1+a-c-d)}\left\{\frac{a}{2e}\sum_{n=0}^\infty \frac{(c)_n(d)_n(\fs+\fs a-b)_n}{(\fs+\fs a)_n (1+a-b)_n n!}\right.\hspace{3cm}\]
\[\hspace{4cm}\left.+\left(1-\frac{a}{2e}\right)\sum_{n=0}^\infty\frac{(c)_n (d)_n (1+\fs a-b)_n}{(1+\fs a)_n (1+a-b)_n n!}\right\}.\]
Identification of the resulting series on the right-hand side as ${}_3F_2(1)$ series then leads to the result stated in (\ref{e21}). The convergence of the hypergeometric series appearing on the left and right-hand sides of (\ref{e21}) requires
$\Re (a-b-c-d)>-\f{3}{2}$ and $\Re (a-c-d)>-1$, respectively; combination of these two conditions yields the stated condition following (\ref{e21}).
This completes the proof of the theorem.\ \ \ \ \ $\Box$
\vspace{0.2cm}

If we put $e=\fs a$ in (\ref{e21}), we recover Exton's result in (\ref{e16}). Thus (\ref{e21}) may be regarded as a generalization of (\ref{e16}), which we hope may prove to be of interest.

Further, if $d=1+a-b$ in (\ref{e21}), we find
\[{}_3F_2\left[\begin{array}{c}a,\ c,\ 1+e\\1+a-c,\ e\end{array};-1\right]\hspace{9cm}\]
\[=\frac{\g(1+a-c)}{\g(1+a)}\left\{\frac{a}{2e} \,{}_2F_1\left[\begin{array}{c}\vspace{0.1cm}

c,\,\fs+\fs a-b\\\fs+\fs a\end{array}\!;1\right]+\left(1-\frac{a}{2e}\right)\,{}_2F_1\left[\begin{array}{c}\vspace{0.1cm}

c,\,1+\fs a-b\\1+\fs a\end{array}\!;1\right]\right\}\]
\bee\label{e34}
=\frac{\g(1+a-c)}{\g(1+a)}\left\{\frac{a}{2e}\,\frac{\g(\fs+\fs a)}{\g(\fs+\fs a-c)}+\left(1-\frac{a}{2e}\right)\,\frac{\g(1+\fs a)}{\g(1+\fs a-c)}\right\}
\ee
provided $\Re(c)<\fs$. [We note that in the evaluation of the ${}_2F_1(1)$ series by Gauss' theorem we need the dummy condition $\Re(b-c)>0$.] This formula has been obtained by different means in \cite[Eq.~(5.10)]{KRR} and \cite[Eq.~(4.7)]{RP}. If $e=\fs a$, (\ref{e34}) reduces to the summation formula given in \cite[p.~245, III.21]{S}; if, in addition, $c=-n$, where $n$ is a non-negative integer, (\ref{e34}) reduces to (III.25) in \cite[p.~245]{S}.
\vspace{0.5cm}

\noindent{\bf Acknowledgement:}\ \ \ Y. S. Kim acknowledges the support of the Wonkwang University Research Fund (2015).

\end{document}